\documentclass{article}

\usepackage[utf8]{inputenc}
\usepackage{hyperref}
\usepackage{amssymb}
\usepackage{amsmath}
\usepackage[english]{babel}
\newcommand{\M}{{\mathcal M}}
\newtheorem{lemma}{Lemma}
\newtheorem{theorem}{Theorem}
\usepackage{tikz}
\usetikzlibrary{shapes,arrows}
\tikzset{cblue/.style={circle, draw, thin,fill=cyan!20, scale=0.5}}
\tikzset{cred/.style={circle, draw, thin, fill=red!20, scale=0.5}}
\tikzset{cgreen/.style={circle, draw, thin, fill=green!20, scale=0.5}}
\tikzset{cblack/.style={circle, draw, thin, fill=black, scale=0.2}}
\date{}
\author{E.Yu.~Lerner, R.E.~Lerner}
\title{Minimal instances with no weakly stable matching for three-sided problem with cyclic incomplete preferences}
\begin{document}
\maketitle
\begin{abstract}
Given $n$ men, $n$ women, and $n$ dogs, each man has an incomplete preference list of women, each woman does an incomplete preference list of dogs, and each dog does an incomplete preference list of men.
We understand a family as a triple consisting of one man, one woman, and one dog such that each of them enters in the preference list of the corresponding agent.
We do a matching as a collection of nonintersecting families (some agents, possibly, remain single).
A matching is said to be nonstable, if one can find a man, a woman, and a dog which do 
not live together currently but each of them would become ``happier'' if they do. Otherwise the matching is said to be stable (a weakly stable matching in 3-DSMI-CYC problem). We give an example of this problem for $n=3$ where no stable matching exists. Moreover, we prove the absence of such an example for $n<3$. 
Such an example was known earlier only for $n=6$ (Biro, McDermid, 2010). The constructed examples also allows one to decrease 
(in two times) the size of the recently constructed analogous example for complete preference lists (Lam, Plaxton, 2019).
\end{abstract}
\section{Introduction}
Assume that there are $n$ men and $n$ women, and each one among them has a preference list of representatives of the opposite sex. A partition into heterogeneous families with no tuple of man and woman, who prefer each other rather than their partners (if they have ones), is called stable matching. The initial case of complete preference lists was studied by D.~Gale and L.S.~Shaple: a stable matching necessarily exists, an $O(n^2)$-hard algorithm for forming it was proposed in~\cite{Gale}. Note that the algorithm is also applicable for the case of incomplete preference lists, but some men and women, possibly, remain single. A certain modification of the Gale--Shaple algorithm (see, for example,~\cite{knuth}) allows one to find, if possible, a matching without single men and women or to prove its absence, otherwise. 

In the case of random complete preference lists (when for each man and each woman the distribution of all permutations of representatives of the opposite sex is independent and uniform) the time necessary for finding a stable matching is $\Theta( n\ln(n))$~\cite{knuth}.
In this case, the mean value of the number of stable matchings also has the asymptote $n\ln(n)$~\cite{pittel1}.

In~\cite{knuth} D.~Knuth states the question whether it is possible to generalize the theory of stable matchings to the case of three genders. The most interesting variant in the $k$-gender case occurs when preferences are cyclic: representatives of the 1st gender rank representatives of the 2nd one, the latter do representatives of the 3rd gender, etc., and each representative of the $k$th gender has a preference list of representatives of the 1st gender (see~\cite[Chapter~5.6]{manlove} for the non-cyclic variants of the $k$-gender case).

A tuple containing exactly one representative of each gender is called a family, and the set of disjoint families is called a matching. A matching is said to be weakly stable, if there is no tuple outside this matching, each member of which would become ``happier'', if they live together. In what follows, for brevity, we use the term ``a stable matching'' instead of the term ``a weakly stable matching''.

Let the number of representatives of each gender equal~$n$. In~\cite{gurvich}, one proves that with complete preference lists a stable matching always exists, provided that $n\leqslant k$ (where $k$ is the number of genders). In~\cite{Eriksson}, Eriksson et al. generalize this result for the case when $k=3$ and $n=k+1=4$. Ibid, one states the conjecture that the problem of finding a stable matching in 3-gender case with complete preference lists (problem~3-DSM-CYC or just 3DSM) has a solution for any~$n$.
Using a satisfiability problem formulation and an extensive computer-assisted search, the authors of~\cite{new} prove the validity of the conjecture stated by Eriksson et al. for~$n=5$. In~\cite{Pittel2}, one proves that with random preference lists the mean value of stable matchings in problem 3DSM grows as $\Omega(n^2 \ln^2(n))$.

The 3DSMI-problem (3-dimensional stable matching with incomplete preference lists) was studied by P.~Bir\'o and E.~McDermid~\cite{Biro}. According to results a solution of 3DSMI does not necessarily exists in contrast to the two-dimensional case; they give an explicit example of problem 3DSMI for $n=6$ with no stable matching. Moreover, they prove that the problem of establishing the solvability of 3DSMI is NP-complete. Ibid, they state the problem of constructing an instance with no weakly stable matching for $n<6$.

Finally, contrary to expectations, the conjecture stated by~Eriksson et al. was recently refuted in~\cite{Lam}. Lam and Paxton associate problem 3DSMI with a certain problem 3DSM, where $n$ is 15 times greater than the initial size; this problem is solvable if and only if so is the initial problem 3DSMI. Therefore, the problem of establishing the solvability of problem 3DSM is NP-complete. The example described in the paper~\cite{Biro} allows one to construct an instance of problem 3DSM for $n=90=6\times15$ with no stable matching.

For this reason, the problem of finding an instance of 3DSMI with no weakly stable matching for 
$n<6$ becomes more actual. 
The construction of such instances for the least possible values of $n$ is the goal of this paper.

First we constructed an instance of 3DSMI problem for $n=4$ and proved the absence of such instances for $n<3$. But after failing to prove the absence of such instances for $n=3$, we have proposed an algorithm for the computer search of all possible instances of 3DSMI problems for $n=3$. Unexpectedly, the algorithm has succeeded in constructing instances of rather simple 3DSMI problems without weakly stable matching for $n=3$.

The rest part of the paper has the following structure.
In Sect.~2, we present the formal definitions of 3DSMI-CYC in terms of the graph theory.
In Sect.~3, we study some properties of graphs of problem~3DSMI-CYC,
prove the absence of counterexamples for $n<3$.
In Sect.~4 we describe various cases of problem 3DSMI for $n=3$ and
consider the result of their computer enumeration. 
We consider several instances and explicitly prove the absence of a stable matching for each of them, 
and describe general properties of all counterexamples.
In Sect.~5, we conclude by mentioning some potential future work.
In Appendix, we describe our example for $n=4$ and prove that in this
case no stable matching exists.

\section{The statement of 3DSMI-CYC in terms of the graph theory}
Let $G$ be some directed graph. Denote the set of its edges by $E$; assume that no edges are multiple. 
Let the vertex set $V (G) = V$ of the graph $G$ be divided into
three subsets, namely, the set of men $M$, women $F$, and dogs $D$. Any edge $(v, v'), v, v' in V$ of this graphis considered to be of one of three types: either $v in M, v' \in F$,
or $v \in F, v' \in D$, or $v \in D, v' \in M$.

Assume that $|M|=|F|=|D|$ (otherwise we supplement the corresponding subgraph with vertices that are not connected with the rest part of the graph). The number $n=|M|=|F|=|D|$ is called the problem \textit{size}. Evidently, the length of all cycles in the graph~$G$ is a multiple of~$3$. Note also that this condition ensures the possibility to divide the vertex set of any digraph~$G$ into 3 subsets $M$, $F$, $D$ so that all its edges are directed as indicated above.
Each edge $(v,v')$, $v,v'\in V$, corresponds to some positive integer $r(v,v')$ which is called the rank of this edge. For fixed $v\in V$, all possible ranks $r(v,v_1),\ldots,r(v,v_k)$ coincide with $\{1,\ldots,k\}$, where $k$ is the outgoing vertex degree~$v$ (if $r(v, v') = 1$, then $v'$ is the best preference for $v$, and so on).

We understand \textit{a three-sided matching} as a subgraph~$H(V)$ of the graph~$G$, where each vertex $v\in V$ has at most one outgoing edge and the following condition is fulfilled: if a vertex $v$ has an outgoing edge, then this edge belongs to a cycle of length~3 in the graph~$H$. Cycles of length~3 in the graph $H$ are called \textit{families}. Evidently, each family, accurate to a cyclic shift, takes the form $(m,f,d)$, where $m\in M$, $f\in F$, and $d\in D$. Note that in what follows, for convenience of denotations of families, we do not fix the order of genders in a family, i.\,e., we treat denotations of families as triples derived from an initial one by a cyclic shift as equivalent.

In what follows, we sometimes use the notion of a family in a wider sense, namely, as any cycle of length~3 in the graph~$G$. However, if some three-sided matching~$H$ is fixed, then we describe other cycles of length~3 explicitly, applying the term ``a family'' only to cycles that enter in a three-sided matching.

The graph G is said to be \textit{trivial} if it has no cycles of length three and \textit{untrivial} otherwise. In the paper we mostly consider graphs to be untrivial.

\textit{A matching} $\M$ is a collection of all families of a three-sided matching $H$. For a vertex $v$, $v\in V$, in the matching $\M$, the rank $R_\M(v)$ is defined as the rank of the edge that goes out of this vertex in the subgraph $H$. If some vertex $v$ in the subgraph $H$ has no outgoing edge, then $R_\M(v)$ is set to $+\infty$.

A triple $(v,v',v'')$ is said to be \textit{blocking} for some matching $\M$, if it is a cycle in the graph~$G$, and
$$
r(v,v')<R_\M(v),\quad r(v',v'')<R_\M(v'),\quad r(v'',v)<R_\M(v'').
$$
A matching $\M$ is said to be \textit{stable}, if no blocking triple exists for it.

Problem~\textit{3DSMI} (3-dimensional stable matching with incomplete preference lists) consists in finding a stable matching for a given graph~$G$. It is well known that it does not necessarily exists. Moreover, the problem of establishing its existence for a given graph~$G$ is NP-complete. It was mentioned in the Introduction, this fact was proved by Biro and McDermid. They have constructed an explicit example of the graph~$G$ of size~6, for which no stable matching exists. Moreover, the question of constructing similar examples for lesser sizes was also stated by the mentioned authors.

\section{The absence of examples of problem 3DSMI with no stable matching for $n<3$}
Let $G$ and $G'$ be two directed graphs defined on one and the same vertex set~$V$ but, generally speaking, having distinct edge sets. Assume that rank functions $r_G$ and $r_{G'}$ are defined on $E$ and $E'$, correspondingly. Let $L\subseteq E\cap E'$. We say that ranking orders $r_G$ and $r_{G'}$ coincide on~$L$, if for any two edges $(v,v'),(v,v'')$ in $L$,
$$
r_G(v,v')<r_G(v,v'')\quad \iff \quad r_{G'}(v,v')<r_{G'}(v,v'').
$$
\begin{lemma}
\label{one}
For any untrivial graph $G$ of problem~3DSMI of size $n$ there exists a graph $G'$ of the same size such that the outgoing degree of each its vertex is nonzero and there is the following correspondence between graphs $G$ and $G'$:\\
1) the set of all possible families of graphs~$G$ and $G'$ coincide;\\
2) the ranking order of all edges that enter in these families also coincide.
\end{lemma}
\textbf{Proof:} Let $v$ be a vertex in the graph~$G$ having no outgoing edges. Then $v$ enters in no family of the graph~$G$. Let us delete this vertex together with all incoming edges in~$v$. Repeating this procedure several times, we get a graph $\widehat{G}$ such that each its vertex has at least one outgoing edge and its set of families coincides with that of the initial graph~$G$. Let the symbol $\widehat V$ stand for the vertex set of the graph~$\widehat G$, denote the set of its edges by $\widehat E$. 
According to untriviality of graph~$G$ the set of families of graphs $G$ and $\widehat G$ is nonempty. In this case, the set $\widehat V$ contains at least one vertex for each gender.

Let us now restore the initial vertices belonging to the set $V\setminus\widehat V$ and for each of them arbitrarily construct at least one edge directed to some vertex in $\widehat V$ that corresponds to a proper gender. Since the incoming degree of restored vertices equals zero, they, as earlier, can enter in no family. Note that $\widehat E\subseteq E$ and, consequently, one can construct a rank function for the obtained graph $G'$ preserving the ranking order of the graph~$G$ on $\widehat E$. The obtained graph~$G'$ with the rank function defined in the indicated way is the desired one.  \qquad  $\square$

Lemma~\ref{one} allows one, when studying problems~3DSMI of size~$n$, to restrict oneself to considering the corresponding graphs~$G$ with nonzero outgoing degrees of all vertices, which we do in what follows.

Let the symbol $G''$ stands for a subgraph of the graph~$G$ consisting of its edges of rank~1. We call $G''$ the \textit{basic subgraph} of the graph~$G$. Since each vertex in the basic subgraph has exactly one outgoing edge, $G''$ represents a collection of cycles, whose lengths are multiples of~3, and trees of edges that lead to these cycles.
\begin{theorem}
\label{two}
Problem~3DSMI of size $n\leqslant 2$ always has a stable matching.
\end{theorem}
\textbf{Proof:} Note that with $n=1$ the assertion of the lemma is trivial. 
In what follows, we restrict ourselves to considering the case of $n=2$. Note also that in this case a nonstable matching can contain only one family.
So let us assume the absence of two family matching for graph~$G$.

The basic subgraph of the graph $G$ contains cycles either of length~3 or of length~6. Let us consider both cases sequentially.
In the first case, the exist vertices $v_0,v_1,v_2$ such that $r( v_0,v_1)=r(v_1,v_2)=r(v_2,v_0)=1$. Therefore, if the family $(v_0,v_1,v_2)$ is a stable matching.

It remains to consider the case when the basic subgraph of the graph~$G$ is a cycle of length~6, i.e., $C=(v_0,v_1,\ldots,v_5)$. Without loss of generality, we assume that the graph~$G$, which represents a counterexample to Theorem~\ref{two}, along with the cycle~$C$ contains the edge $(v_2,v_0)$ of rank~2. Then the only possible blocking triple to the matching of one family $(v_0,v_1,v_2)$ is $(v_2,v_3,v_4)$. Consequently, the graph $G$ also contains the edge $(v_4,v_2)$. But then the only possible blocking triple for the matching consisting of one family $(v_2,v_3,v_4)$ is $(v_4,v_5,v_0)$. In turn, the graph~$G$ that consists of only a basic cycle $C$ and edges $(v_0,v_4)$, $(v_4,v_2)$, $(v_2,v_0)$ of rank~2 has a stable matching consisting of one family $(v_0,v_4,v_2)$. Therefore, the graph~$G$, along with the cycle~$C$, contains at least 4 edges. Consequently, the graph $G$ of size $n=2$ has a matching of two families, and it is stable by definition. \qquad  $\square$

\section{The examples of graphs $G$ of size $n=3$ with no stable matching}
In this section, we consider the case of $n=3$. Let us first classify all graphs of the problem of this size; this will facilitate their computer search.

If the basic subgraph of the graph $G$ contains cycles of length~3, then there exist vertices $v_0,v_1,v_2$ such that $r( v_0,v_1)=r(v_1,v_2)=r(v_2,v_0)=1$. Therefore, if the family $(v_0,v_1,v_2)$ enters in a matching, then these vertices can enter in no blocking triple. But then we again get the case of $n=2$ which is mentioned in the statement of Theorem~\ref{two}.

\begin{figure}[h]
\begin{center}
{
\begin{tikzpicture}[->,>=stealth',shorten >=1pt,auto,node distance=3cm,thick,main node/.style={rectangle,fill=blue!20,draw,font=\sffamily\Large\bfseries}]
        \node[cblack] (0) at ( 0:1) {};
        \node[cblack] (1) at ( 40:1) {};
        \node[cblack] (2) at ( 80:1) {};
        \node[cblack] (3) at ( 120:1) {};
        \node[cblack] (4) at ( 160:1) {};
        \node[cblack] (5) at ( 200:1) {};
        \node[cblack] (6) at ( 240:1) {};
        \node[cblack] (7) at ( 280:1) {};
        \node[cblack] (8) at ( 320:1) {};
        \path[every node/.style={font=\sffamily\small}]
        (0) edge  []  node [] {} (1)
        (1) edge  []  node [] {} (2)
        (2) edge  []  node [] {} (3)
        (3) edge  []  node [] {} (4)
        (4) edge  []  node [] {} (5)
        (5) edge  []  node [] {} (6)
        (6) edge  []  node [] {} (7)
        (7) edge  []  node [] {} (8)
        (8) edge  []  node [] {} (0)
        ;
\end{tikzpicture}
}
{
\begin{tikzpicture}[->,>=stealth',shorten >=1pt,auto,node distance=3cm,thick,main node/.style={rectangle,fill=blue!20,draw,font=\sffamily\Large\bfseries}]
        \node[cblack] (0) at ( 30:1) {};
        \node[cblack] (1) at ( 90:1) {};
        \node[cblack] (2) at ( 150:1) {};
        \node[cblack] (3) at ( 210:1) {};
        \node[cblack] (4) at ( 270:1) {};
        \node[cblack] (5) at ( 330:1) {};
        \node[cblack] (6) at ( 90:1.5) {};
        \node[cblack] (7) at ( 90:2) {};
        \node[cblack] (8) at ( 90:2.5) {};
        \path[every node/.style={font=\sffamily\small}]
        (0) edge  []  node [] {} (1)
        (1) edge  []  node [] {} (2)
        (2) edge  []  node [] {} (3)
        (3) edge  []  node [] {} (4)
        (4) edge  []  node [] {} (5)
        (5) edge  []  node [] {} (0)
        (6) edge  []  node [] {} (1)
        (7) edge  []  node [] {} (6)
        (8) edge  []  node [] {} (7)
        ;
\end{tikzpicture}
}
{
    \begin{tikzpicture}[->,>=stealth',shorten >=1pt,auto,node distance=3cm,thick,main node/.style={rectangle,fill=blue!20,draw,font=\sffamily\Large\bfseries}]
        \node[cblack] (0) at ( 30:1) {};
        \node[cblack] (1) at ( 90:1) {};
        \node[cblack] (2) at ( 150:1) {};
        \node[cblack] (3) at ( 210:1) {};
        \node[cblack] (4) at ( 270:1) {};
        \node[cblack] (5) at ( 330:1) {};
        \node[cblack] (6) at ( 30:1.5) {};
        \node[cblack] (7) at ( 30:2) {};
        \node[cblack] (8) at ( 90:1.5) {};
        \path[every node/.style={font=\sffamily\small}]
        (0) edge  []  node [] {} (1)
        (1) edge  []  node [] {} (2)
        (2) edge  []  node [] {} (3)
        (3) edge  []  node [] {} (4)
        (4) edge  []  node [] {} (5)
        (5) edge  []  node [] {} (0)
        (6) edge  []  node [] {} (0)
        (7) edge  []  node [] {} (6)
        (8) edge  []  node [] {} (1)
        ;
    \end{tikzpicture}
}
{
    \begin{tikzpicture}[->,>=stealth',shorten >=1pt,auto,node distance=3cm,thick,main node/.style={rectangle,fill=blue!20,draw,font=\sffamily\Large\bfseries}]
        \node[cblack] (0) at ( 30:1) {};
        \node[cblack] (1) at ( 90:1) {};
        \node[cblack] (2) at ( 150:1) {};
        \node[cblack] (3) at ( 210:1) {};
        \node[cblack] (4) at ( 270:1) {};
        \node[cblack] (5) at ( 330:1) {};
        \node[cblack] (6) at ( 330:1.5) {};
        \node[cblack] (7) at ( 90:2) {};
        \node[cblack] (8) at ( 90:1.5) {};

        \path[every node/.style={font=\sffamily\small}]
        (0) edge  []  node [] {} (1)
        (1) edge  []  node [] {} (2)
        (2) edge  []  node [] {} (3)
        (3) edge  []  node [] {} (4)
        (4) edge  []  node [] {} (5)
        (5) edge  []  node [] {} (0)

        (6) edge  []  node [] {} (5)
        (7) edge  []  node [] {} (8)
        (8) edge  []  node [] {} (1)
        ;
    \end{tikzpicture}
}
{
    \begin{tikzpicture}[->,>=stealth',shorten >=1pt,auto,node distance=3cm,thick,main node/.style={rectangle,fill=blue!20,draw,font=\sffamily\Large\bfseries}]
        \node[cblack] (0) at ( 0:1) {};
        \node[cblack] (1) at ( 60:1) {};
        \node[cblack] (2) at ( 120:1) {};
        \node[cblack] (3) at ( 180:1) {};
        \node[cblack] (4) at ( 240:1) {};
        \node[cblack] (5) at ( 300:1) {};
        \node[cblack] (6) at ( 300:1.5) {};
        \node[cblack] (7) at ( 180:1.5) {};
        \node[cblack] (8) at ( 60:1.5) {};
        \path[every node/.style={font=\sffamily\small}]
        (0) edge  []  node [] {} (1)
        (1) edge  []  node [] {} (2)
        (2) edge  []  node [] {} (3)
        (3) edge  []  node [] {} (4)
        (4) edge  []  node [] {} (5)
        (5) edge  []  node [] {} (0)
        (6) edge  []  node [] {} (5)
        (7) edge  []  node [] {} (3)
        (8) edge  []  node [] {} (1)
        ;
    \end{tikzpicture}
}
{
    \begin{tikzpicture}[->,>=stealth',shorten >=1pt,auto,node distance=3cm,thick,main node/.style={rectangle,fill=blue!20,draw,font=\sffamily\Large\bfseries}]
        \node[cblack] (0) at ( 0:1) {};
        \node[cblack] (1) at ( 60:1) {};
        \node[cblack] (2) at ( 120:1) {};
        \node[cblack] (3) at ( 180:1) {};
        \node[cblack] (4) at ( 240:1) {};
        \node[cblack] (5) at ( 300:1) {};
        \node[cblack] (6) at ( 300:1.5) {};
        \node[cblack] (7) at ( 240:1.5) {};
        \node[cblack] (8) at ( 180:1.5) {};
        \path[every node/.style={font=\sffamily\small}]
        (0) edge  []  node [] {} (1)
        (1) edge  []  node [] {} (2)
        (2) edge  []  node [] {} (3)
        (3) edge  []  node [] {} (4)
        (4) edge  []  node [] {} (5)
        (5) edge  []  node [] {} (0)
        (6) edge  []  node [] {} (5)
        (7) edge  []  node [] {} (4)
        (8) edge  []  node [] {} (3)
        ;
    \end{tikzpicture}
}
\caption{6 variants of the basic subgraph of the graph $G$.}
\label{G''}
\end{center}
\end{figure}
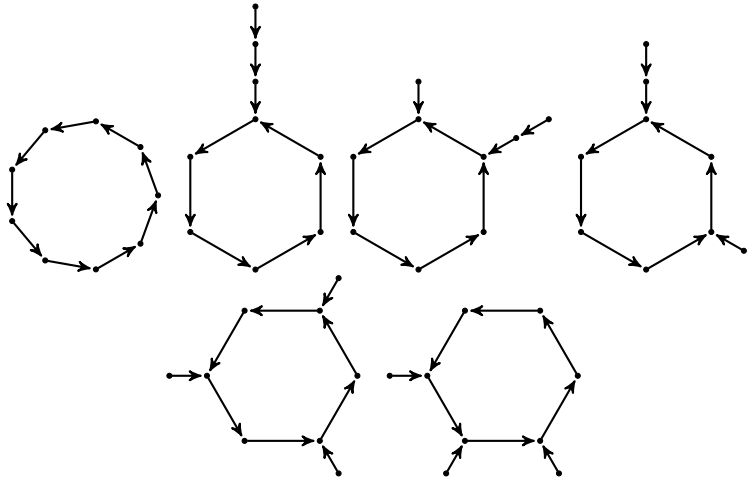

Therefore, the basic subgraph of the graph $G$ represents either a cycle of length~9, or a cycle of length~6 with three edges that lead to this cycle. Altogether, accurate to the cyclic symmetry, there are 6 such subgraphs; they are shown in 
Fig.~\ref{G''}

Each of 9 vertices of these subgraphs in the graph~$G$ can have outgoing edges that lead to two remaining vertices of the corresponding gender (here we understand remaining vertices as those that differ from the vertex, to which the edge of the basic graph~$G''$ is already directed).
Generally speaking, the total number of possible cases is 5, namely, \\
1) the considered vertex has no more outgoing edges;\\
2)-3) the considered vertex has one more outgoing edge that leads to some vertex among two ones;\\
4)-5) the considered vertex has two edges, their ranks are equal to~2 and~3, we can associate ranks with these edges in two ways. \\
Therefore, it suffices to consider $6\times 5^9$ problems 3DSMI.

Evidently, for each of these problems there exist at most 27 families (27 blocking triples). The number of possible three-sided matching, as one can easily calculate, also is not so large. Namely, there exist at most 27 matchings consisting of one family. In addition, there exist at most 108 matchings consisting of two families, namely,
there are 27 ways to form a triple consisting of repre\-sen\-ta\-tives of genders that enter in no matching, and 4 ways to choose partners among two women and two dogs entering the matching for a fixed man that also enters this matching.
Finally, there is at least 36 three-sided matchings of 3 triples, 
as there are $3! \times 3!$ ways to distribute women and dogs among three man-indexed triplets.
Therefore, the total amount of three-sided matchings does not exceed 36+108+27=171.
For each of them we need to find the first triple among 27 potential blocking ones that really is blocking.

Therefore, the total amount of considered cases does not exceed $6\times5^9\times 171\times 27\approx 54\times 10^9$. For generating these cases, we have written a program in Python. 
See the version of this program that calculates the number of counterexamples for each of basic graphs shown in Fig.~\ref{G''} at \url{https://github.com/reginalerner/3dsm/}.

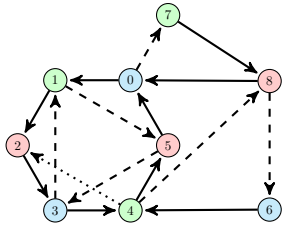
\begin{figure}[h]
\begin{center}
\begin{tikzpicture}[->,>=stealth',shorten >=1pt,auto,node distance=3cm,thick,main node/.style={rectangle,fill=blue!20,draw,font=\sffamily\Large\bfseries}]
        \node[cred] (5) at ( 0:1) {5};
        \node[cblue] (0) at ( 60:1) {0};
        \node[cgreen] (1) at ( 120:1) {1};
        \node[cred] (2) at ( 180:1) {2};
        \node[cblue] (3) at ( 240:1) {3};
        \node[cgreen] (4) at ( 300:1) {4};
        \node[cblue] (6) at ( 340:2.5) {6};
        \node[cgreen] (7) at ( 60:2) {7};
        \node[cred] (8) at ( 20:2.5) {8};

        \path[every node/.style={font=\sffamily\small}]
        (0) edge  []  node [] {} (1)
        (1) edge  []  node [] {} (2)
        (2) edge  []  node [] {} (3)
        (3) edge  []  node [] {} (4)
        (4) edge  []  node [] {} (5)
        (5) edge  []  node [] {} (0)

        (6) edge  []  node [] {} (4)
        (7) edge  []  node [] {} (8)
        (8) edge  []  node [] {} (0)
        
        (4) edge  [dashed]  node [] {} (8)
        (8) edge  [dashed]  node [] {} (6)
        (0) edge  [dashed] node [] {} (7)
        (4) edge  [dotted]  node [] {} (2)    
        (1) edge  [dashed]  node [] {} (5)
        (5) edge  [dashed]  node [] {} (3)
        (3) edge  [dashed] node [] {} (1)
        ;
\end{tikzpicture}
\caption{The graph of problem 3DSMI of size~3 with no stable matching consisting of 16 edges. 
The rank of all edges indicated by solid bold lines equals~1. The dashed lines represent the edges with the rank~2. The rank of the ``dotted'' edge equals~3. 
}
\label{example}
\end{center}
\end{figure}

For the first basic graph shown in Fig.~\ref{G''} (a cycle formed by 9 vertices), no such graph for problem 3SDMI with no stable matching was obtained. As is mentioned in the Introduction, we even did not expect to find such instances for $n=3$. To our surprise, the computer search has found such counterexamples for each of the rest basic graphs. One of them is shown in Fig.~\ref{example}.
For convenience, we enumerate vertices of the graph by numbers~$v$, $v=0,1,\ldots,8$. The value $v\bmod 3$ defines the gender that corresponds to the vertex~$v$ .

\begin{theorem}
\label{main}
There is no stable matching in 3DSMI problem for the the graph on Fig.~\ref{example}.
\end{theorem}
\textbf{Proof:}
Fig.~\ref{example} evidently demonstrates that each possible cycle of length 3 takes one of the following forms: $(0,1,5)$, $(0,7,8)$, $(1,2,3)$, $(1,5,3)$, $(2,3,4)$, $(3,4,5)$, and $(4,8,6)$.
These cycles form families, while collections of disjoint families form matchings~$\M$ in the problem.

Evidently, if one can add a cycle $(v,v',v'')$ to a matching $\M$ (vertices $v,v',v''$ do not enter in $\M$), then $\M$ is unstable, i.e., the triple $(v,v',v'')$ is blocking for~$\M$. Therefore, candidates for stable matchings should be supplemented with possible cycles. We call such matchings \textit{uncompletable} and consider only such ones.

The union of vertices of three cycles that are listed above does not coincide with the set of all vertices of the graph shown in Fig.~\ref{example}. On the other hand, by using the direct search method we can prove that any set consisting of one triple is completable. Therefore, each uncompletable matching consists of two families. 
Below we give their complete list together with blocking triples:\\
1) $\{(0,1,5),(2,3,4)\}$, the blocking triple is $(4,8,6)$;\\
2) $\{(0,1,5),(4,8,6)\}$, the blocking triple is $(1,2,3)$;\\
3) $\{(0,7,8),(1,2,3)\}$, the blocking triple is $(3,4,5)$;\\
4) $\{(0,7,8),(1,5,3)\}$, the blocking triple is $(2,3,4)$;\\
5) $\{(0,7,8),(2,3,4)\}$, the blocking triple is $(0,1,5)$;\\
6) $\{(0,7,8),(3,4,5)\}$, the blocking triple is $(0,1,5)$;\\
7) $\{(1,2,3),(4,8,6)\}$, the blocking triple is $(0,7,8)$ or $(3,4,5)$;\\
8) $\{(1,5,3),(4,8,6)\}$, the blocking triple is $(0,7,8)$. \qquad  $\square$

%%%%%%%%%%%%%%%%%%%
One can easily give other examples of graphs with the same set of cycles, uncompletable matchings,
and blocking triples. In particular, this property is characteristic for the graph that differs from that shown in Fig.~\ref{example} by the presence of the additional edge $(7,2)$ of rank~2 or the additional edge $(6,1)$ of rank~2, or both of these edges.
%%%%%%%%%%%

\begin{figure}[h]
\begin{center}
\begin{tikzpicture}[->,>=stealth',shorten >=1pt,auto,node distance=3cm,thick,main node/.style={rectangle,fill=blue!20,draw,font=\sffamily\Large\bfseries}]
        \node[cred] (5) at ( 0:1.2) {5};
        \node[cblue] (0) at ( 60:1.2) {0};
        \node[cgreen] (1) at ( 120:1.2) {1};
        \node[cred] (2) at ( 180:1.2) {2};
        \node[cblue] (3) at ( 240:1.2) {3};
        \node[cgreen] (4) at ( 300:1.2) {4};
        \node[cblue] (6) at ( 150:1.75) {6};
        \node[cgreen] (7) at ( 210:1.75) {7};
        \node[cred] (8) at ( 90:1.75) {8};
        \path[every node/.style={font=\sffamily\small}]
        (0) edge  []  node [] {} (1)
        (1) edge  []  node [] {} (2)
        (2) edge  []  node [] {} (3)
        (3) edge  []  node [] {} (4)
        (4) edge  []  node [] {} (5)
        (5) edge  []  node [] {} (0)
        (6) edge  []  node [] {} (1)
        (7) edge  []  node [] {} (2)
        (8) edge  []  node [] {} (0)
        
        (1) edge  [dashed]  node [] {} (8)
        (2) edge  [dashed]  node [] {} (6)
        (0) edge  [dashed]  node [] {} (4)
        (3) edge  [dashed] node [] {} (7)
        (0) edge  [dotted]  node [] {} (7)    
        (5) edge  [dashed]  node [] {} (3)
        (7) edge  [dashed] node [] {} (5)
        ;
\end{tikzpicture}
\caption{One more 3D graph of problem 3DSMI with no stable matching consisting of 16 edges. 
Denotations are the same as in Fig.~\ref{example}.
}
\label{example1}
\end{center}
\end{figure}
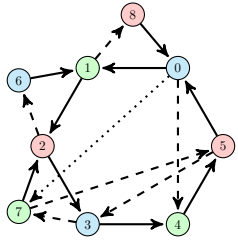

Moreover, one can find other graphs consisting of 16 edges that have no stable matching. One of them is shown in Fig.~\ref{example1} (any other graph with this property differs from the indicated one only in the fact that ranks of edges $(0,4)$ and $(0,7)$ have interchanged).
Note that the graph shown in Fig.~\ref{example1} is similar to that in our example for $n=4$ (see, for comparison, Fig.~\ref{4examp} in Appendix).
These graphs define the following families of forming matchings in problem 3DSMI: $(0,1,8)$, $(0,4,5)$, $(0,7,5)$, $(1,2,6)$, $(2,3,7)$, $(3,4,5)$, and $(3,7,5)$.

The list of matchings with blocking triples looks as follows:\\
1) $\{(0,1,8),(2,3,7)\}$, the blocking triple is $(3,4,5)$;\\
2) $\{(0,1,8),(3,4,5)\}$, the blocking triple is $(1,2,6)$;\\
3) $\{(0,1,8),(3,7,5)\}$, the blocking triple is $(1,2,6)$;\\
4) $\{(0,4,5),(1,2,6)\}$, the blocking triple is $(2,3,7)$;\\
5) $\{(0,4,5),(2,3,7)\}$, the blocking triple is $(0,1,8)$;\\
6) $\{(0,7,5),(1,2,6)\}$, the blocking triple is $(2,3,7)$;\\
7) $\{(1,2,6),(3,4,5)\}$, the blocking triple is $(0,7,5)$.\\
8) $\{(1,2,6),(3,7,5)\}$, the blocking triple is $(0,4,5)$.

Since the counterexamples considered above are diverse, they have some common properties.
We are going to describe them in a future paper.

\section{Concluding remarks}
In this paper, we study the problem stated by Biro and McDermid in~\cite{Biro},
namely, we seek for instances with no weakly stable matching for 3DSM-CYC with $n<6$.
In particular, we find the minimal value of $n$, with which such instances exist,
and describe some of them.

The idea of this study is due to the work of Lam and Paxton~\cite{Lam}, who give an example of problem 3DSM-CYC for $n=90$ with no stable matching. This example is based on an analogous example proposed by Biro and McDermid for problem~3DSMI-CYC with $n=6$. Our example constructed for problem~3DSMI with $n=3$ allows one to make the size of an example for problem~3DSM with no stable matching as low as $n=45$. According to results obtained in Sect.~3, the further decrease of $n$ for problem~3DSMI is impossible. However, it seems possible to find problem 3DSM with no stable matching with $n<45$ using some other methods.

Actually, Lam and Paxton studied not only 3-DSM-CYC, but also its $k$-gender analog, $k$-DSM-CYC, for arbitrary $k\geqslant 3$. First they have represented problem~3-DSMI-CYC as a particular case of $k$-DSMI-CYC with $n^2$ representatives of each gender. Then by the reduction from $k$-DSMI-CYC they have proved that $k$-DSM-CYC is NP-complete.

Note that some development of ideas proposed in the paper~\cite{Lam} allows one to rather easily construct a counterexample of size~$n=5$ 
for $k$-DSMI-CYC, $k>3$, basing on the graph shown in Fig.~\ref{example} via subdivision of outcoming edges of woman vertexes. Any of subdivided edges is converted to the chain with $k-3$ vertexes inside, one for each new gender. A $k$-gender family should contain the new vertexes from subdivided edge, so there is a biunique correspondence between new $k$-gender families and old 3-gender ones. If no stable matching exists for 3-gender families, then neither one exists for the new $k$-gender graph. 

Therefore, for any $k>2$, we have constructed an instance of problem $k$-DSMI-CYC with $n=5$ with no stable matching (where lists of preferences of two women, two men, and two dogs that are not shown in Fig.~\ref{example} can be arbitrary).
The question about the existence of such counterexamples for $n=3$ and $n=4$ still remains open.

We hope that this work can be useful in studying other questions related to other aspects of the generalization of the theory of stable matchings to the $k$-dimensional case, $k>2$. In our opinion, this study is far from completion.

\section*{Appendix. An example of a graph $G$ of size $n=4$ with no stable matching}
\begin{figure}[h]
\begin{center}
        \begin{tikzpicture}[->,>=stealth',shorten >=1pt,auto,node distance=3cm,thick,main node/.style={rectangle,fill=blue!20,draw,font=\sffamily\Large\bfseries}]
        \node[cblue] (2) at ( 0:2) {$v_2$};
        \node[cred] (1) at ( 40:2) {$v_1$};
        \node[cgreen] (0) at ( 80:2) {$v_0$};
        \node[cblue] (8) at ( 120:2) {$v_8$};
        \node[cred] (7) at ( 160:2) {$v_7$};
        \node[cgreen] (6) at ( 200:2) {$v_6$};
        \node[cblue] (5) at ( 240:2) {$v_5$};
        \node[cred] (4) at ( 280:2) {$v_4$};
        \node[cgreen] (3) at ( 320:2) {$v_3$};

        \node[cblue] (w2) at ( 60:3 ) {$w_2$};
        \node[cgreen] (w3) at ( 20:3 ) {$w_3$};
        \node[cred] (w4) at ( 340:3 ) {$w_4$};

        \path[every node/.style={font=\sffamily\small}]
        (0) edge  []  node [] {} (1)
        (1) edge  []  node [] {} (2)
        (2) edge  []  node [] {} (3)
        (3) edge  []  node [] {} (4)
        (4) edge  []  node [] {} (5)
        (5) edge  []  node [] {} (6)
        (6) edge  []  node [] {} (7)
        (7) edge  []  node [] {} (8)
        (8) edge  []  node [] {} (0)

        (w2) edge  []  node [] {} (0)
        (w3) edge  []  node [] {} (1)
        (w4) edge  []  node [right] {} (2)

        (1) edge  [dashed]  node [right] {} (w2)
        (2) edge  [dashed]  node [right] {} (w3)
        (3) edge  [dashed]  node [right] {} (w4)

        (0) edge  [dashed]  node [right] {} (4)
        (1) edge  [dotted]  node [right] {} (5)
        (2) edge  [dotted]  node [right] {} (6)
        (3) edge  [dotted]  node [right] {} (7)
        (4) edge  [dashed]  node [right] {} (8)
        (5) edge  [dashed]  node [right] {} (0)
        (6) edge  [dashed]  node [right] {} (1)
        (7) edge  [dashed]  node [right] {} (2)
        (8) edge  [dashed]  node [right] {} (3)
        ;
        \end{tikzpicture}
\caption{An example of a graph with no stable matching for $n=4$}
\label{4examp}
\end{center}
\end{figure}
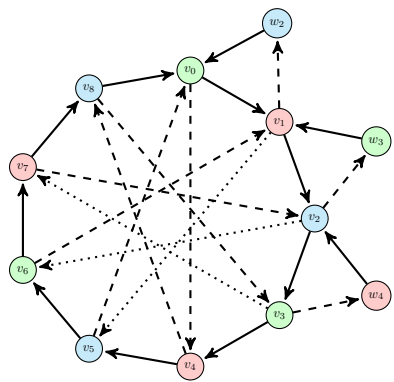
Consider the graph $G$ shown in Fig.~\ref{4examp}. Here bold lines indicate edges of rank~1, dashed ones do edges of rank~2, and dotted ones do edges of rank~3. Therefore, the vertex set of the graph is $$V=\{v_0,\ldots,v_8,w_2,w_3,w_4\},$$ the remainder $i\bmod 3$ of the division of the vertex index~$i$ by 3 defines the gender that corresponds to this vertex.
Edges of the graph~$G$ take one of the following forms:\\
1) $(v_i,v_{(i+1)\bmod 9})$, $i=0,\ldots,8$;\\
2) $(w_i,v_{i-2})$, $i=2,3,4$;\\
3) $(v_i,w_{i+1})$, $i=1,2,3$;\\
4) $(v_i,v_{(i+4)\bmod 9})$, $i=0,\ldots,8$.\\
The rank of edges of the 1st and 2nd kinds equals~1; the rank of edges of the 3d kind equals~2; the rank of edges of the 4th kind equals~2 for $i=0,4,5,\ldots,8$ and~3 for $i=1,2,3$.
In what follows in this section, we consider indices of all vertices modulo~9, for brevity of notations, we omit the symbol $\bmod\,9$ in subscripts.

\begin{theorem}
There is no stable matching in 3DSMI problem for the the graph on Fig.~\ref{4examp}.
\end{theorem}
\textbf{Proof:}
Analogously to the proof of Theorem~\ref{main}, we consider only uncompletable matchings.

One can easily see that all possible families of the graph~$G$ take one of the following forms: $$(v_i,v_{i+1},v_{i+5}), i=0,\ldots,8,\quad \text{or}\quad (v_i,v_{i+1},w_{i+2}), i=0,1,2.$$
All uncompletable matchings for the graph~$G$ consist of two or three such families.

In particular, uncompletable matchings of two families take the form
$$
\{(v_i,v_{i+1},v_{i+5}),(v_{i+2},v_{i+3},v_{i+7})\},\quad  i=0,\ldots,8,
$$
where one or two vertices among $v_{i+5},v_{i+7}$ can be replaced with those $w_{i+2}$ and $w_{i+4}$, correspondingly, (certainly, this replacement is possible, only if a vertex $w$ with the corresponding index exists).
In any case, for such a matching $\M$ the blocking triple takes the form $(v_{i+3},v_{i+4},v_{i+8})$. Really, by definition, $$R_\M(v_{i+4})=R_\M(v_{i+8})=\infty,\quad \text{and}\quad R_\M(v_{i+3})>1.$$

Uncompletable matchings of three families for the graph~$G$ take the form:
$$
\M_1(i)=\{(v_i,v_{i+1},v_{i+5}),(v_{i+2},v_{i+6},v_{i+7}),(v_{i+3},v_{i+4},v_{i+8})\},\quad  i=0,\ldots,8,
$$
or
$$
\M_2(i)=\{(v_i,v_{i+1},w_{i+2}),(v_{i+2},v_{i+6},v_{i+7}),(v_{i+3},v_{i+4},v_{i+8})\},\quad  i=0,1,2.
$$
Note that $\M_2(i)\equiv\M_3((i+3)\bmod 9)\equiv\M_4((i+6)\bmod 9)$, where
$$
\M_3(i)=\{(v_i,v_{i+1},v_{i+5}),(w_{i+8},v_{i+6},v_{i+7}),(v_{i+3},v_{i+4},v_{i+8})\},\quad  i=3,4,5;
$$
$$
\M_4(i)=\{(v_i,v_{i+1},v_{i+5}),(v_{i+2},v_{i+6},v_{i+7}),(v_{i+3},v_{i+4},w_{i+5})\},\quad  i=6,7,8.
$$
In addition, $\M_1(i)\equiv\M_1((i+3)\bmod 9)\equiv\M_1((i+6)\bmod 9)$. Therefore, there exist, as a total, 6 uncompletable matchings of three families:
$$\M_1(0),\M_2(0),\M_1(1),\M_2(1),\M_1(2),\M_2(2).$$

For matchings $\M_1(0),\M_2(0)$ the blocking triple is $(v_1,v_2,w_3)$: $R_{\M_1(0)}(v_1)=3>1$, $R_{\M_2(0)}(v_1)=2>1$,
$R_{\M_1(0)}(v_2)=R_{\M_2(0)}(v_2)=3>2$.
Analogously, for matchings $\M_1(1),\M_2(1)$ the blocking triple is $(v_2,v_3,w_4)$; for matchings $\M_1(2),\M_2(2)$ the blocking triple is $(v_0,v_1,w_2)$.\qquad  $\square$
\end{document}